\documentclass{article}
\usepackage{amsthm}
\usepackage{amsmath}
\newtheorem{thm}{Theorem}[section] 
\newtheorem{lemma}[thm]{Lemma} 
 
\newtheorem{cor}[thm]{Corollary}

\theoremstyle{definition}
\newtheorem{defin}[thm]{Definition}
 \theoremstyle{remark}
\newtheorem{rems}[thm]{Remarks}

\usepackage{amssymb}
\usepackage{hyperref}
\usepackage{color}
\usepackage{calrsfs}            % use a different \mathcal font

\newcommand{\FF}{\mathbb F} 
\newcommand{\UU}{\mathbb U}
\newcommand{\cB}{\mathcal{B}}
\newcommand{\cD}{\mathcal{D}}
\newcommand{\cL}{\mathcal{L}}
\newcommand{\cM}{\mathcal{M}}
\newcommand{\cH}{\mathcal{H}}
\newcommand{\Aut}{\operatorname{Aut}}
\newcommand{\PG}{\mathrm{PG}}

\newcommand{\set}[2]{\{{#1}\,\vert\,{#2}\}}

\makeatletter
\newcommand{\subjclass}[2][1991]{%
  \let\@oldtitle\@title%
  \gdef\@title{\@oldtitle\footnotetext{#1 \emph{Mathematics subject classification.} #2}}%
}
\newcommand{\keywords}[1]{%
  \let\@@oldtitle\@title%
  \gdef\@title{\@@oldtitle\footnotetext{\emph{Key words and phrases.} #1.}}%
}
\makeatother
\begin{document}

\title{Finite subunitals of the Hermitian unitals}

\author{Theo Grundh\"ofer, Markus J. Stroppel, Hendrik {Van Maldeghem}}

\date{}



\maketitle

\begin{abstract}
  \noindent%
  Every finite subunital of any generalized hermitian unital is itself
  a hermitian unital; the embedding is given by an embedding of
  quadratic field extensions.  In particular, a generalized hermitian
  unital with a finite subunital is a hermitian one (i.e., it
  originates from a separable field extension).
\end{abstract}

\enlargethispage{5mm}%
\noindent%
A hermitian unital in a pappian projective plane consists of the
absolute points of a unitary polarity of that plane, with blocks
induced by secant lines (see Section \ref{herm}).  The finite
hermitian unitals of order $q$ are the classical examples of
$2$-$(q^3+1, q+1, 1)$-designs, cp.~\cite[II.8, pp.\,57--63 and
p.\,246]{MR0333959}, or~\cite[Ch.\,2]{MR2440325}. %
We consider generalized hermitian unitals $\cH(C|R)$ where $C|R$ is
any quadratic extension of fields; separable extensions $C|R$ yield
the hermitian unitals, inseparable extensions give certain projections
of quadrics.

If one has an embedding of field extensions from $E|F$ into $C|R$ (in
the sense of~\ref{def:standard}) then clearly the unital $\cH(E|F)$ is
embedded into $\cH(C|R)$. %
Our Main Theorem asserts, in the converse direction, that every
embedding of a finite unital into one of those generalized hermitian
unitals comes from an embedding of field extensions.  This also
implies that no finite subunitals exist in generalized hermitian
unitals that are not hermitian unitals.

If $\UU$ is a unital of order~$q$ embedded in the hermitian unital
$\cH(\FF_{r^2}|\FF_r)$ of order~$r$ then our result says that~$\UU$ is
the hermitian unital of order~$q$, and there is an odd integer~$e$
such that $r=q^e$.

We remark that an analogous result holds for generalized polygons
(where the Moufang property singles out the classical examples),
see~\cite[5.2.2]{MR1725957}: any (thick) subpolygon of a Moufang
polygon is Moufang, as well.

\section{Generalized hermitian unitals}\label{herm}

Let $C|R$ be any quadratic (possibly inseparable) extension of fields;
so $C=R+\varepsilon R$, with $\varepsilon\in C\smallsetminus R$.  There
exist $t,d\in R$ such that $\varepsilon^2-t\varepsilon +d=0$, since
$\varepsilon^2\in R+\varepsilon R$.  The mapping
\[
\sigma\colon C\rightarrow C \colon %
x+\varepsilon y\mapsto (x+ty)-\varepsilon y \quad\text{for } x,y\in R
\]
is a field automorphism which generates $\Aut_R C$: if $C|R$ is
separable, then $\sigma$ has order $2$; if $C|R$ is inseparable, then
$\sigma$ is the identity.

Now we introduce our geometric objects. We consider the pappian
projective plane $\PG(2,C)$ arising from the $3$-dimensional vector
space $C^3$ over $C$, and we use homogeneous coordinates
$[X,Y,Z] := (X,Y,Z)C$ for the points of $\PG(2,C)$. %

\goodbreak%
\begin{defin}
  The \emph{generalized hermitian unital} $\cH(C|R)$ is the incidence
  structure %
  $(H, \cM)$ with the point set %
  $H:= \set{ [X,Y,Z] }{ X^\sigma Y+Z^\sigma Z\in\varepsilon R}$, and
  the set $\cM$ of \emph{blocks} consists of the intersections of $U$
  with secant lines, i.e.\ lines of~$\PG(2,C)$ containing more than
  one point of $H$.
\end{defin}

\begin{rems}\label{identify}
  The name ``generalized hermitian unital'' is motivated by the
  following observations.
  
  If\/ $C|R$ is separable, then $\cH(C|R)$ is the hermitian unital
  arising from the skew-hermitian form %
  $h\colon C^3\times C^3\to C$ defined by
  \[
    h\left((X,Y,Z),(X',Y',Z')\right) = %
    \varepsilon^\sigma X^\sigma Y'-\varepsilon Y^\sigma X'
    +(\varepsilon^\sigma-\varepsilon)Z^\sigma Z' \,,
  \]
  see \cite[2.2]{MR4008656}. %
  So the elements of~$H$ are the absolute points of a polarity of
  $\PG(2,C)$, in this case.
  
  If\/ $C|R$ is inseparable, then $\cH(C|R)$ is the projection of an
  ordinary quadric $Q$ in some projective space of dimension at least
  $3$ from a subspace of codimension~$1$ in the nucleus of $Q$, see
  \cite[2.2]{MR4008656}.  In this case, there is no polarity of
  $\PG(2,C)$ such that the elements of~$H$ are absolute points.
  
  In any case, we have for each point~$p$ of $\cH(C|R)$ a unique line of
  $\PG(2,C)$ passing through~$p$ and containing no other point of
  $\cH(C|R)$,
   see \cite[2.3]{MR4008656}. We will refer to this line as  the
  tangent line to~$\cH(C|R)$ at~$p$ in~$\PG(2,C)$.

  The unital $\cH(C|R)$ admits all conceivable translations (i.e.,
  automorphisms that fix a point and every block through it), and
  these are induced by elations of~$\PG(2,C)$, see
  \cite[2.13]{MR4008656}.

  Each block of~$\cH(C|R)$ is a Baer subline in~$\PG(2,C)$, %
  see \cite[2.8]{MR4008656}, cp.~\cite{MR4224387}. %
\end{rems}

By an \emph{O'Nan configuration}, we mean four blocks intersecting in
six points of the unital (i.e., a $(6_2\,4_3)$ configuration). %
Naming this configuration in honor of O'Nan~\cite{MR0295934} is
customary in the context of unitals, see~\cite[p.\,87]{MR2440325}; the
configuration is named after Veblen and Young in the axiomatics of
projective spaces, or after Pasch in the context of ordered
(Euclidean) geometry.

\bigskip%
We will use the following properties of generalized hermitian
unitals: %

%\enlargethispage{20mm}%
\begin{lemma}\label{properties}%
Let $C|R$ be a quadratic field extension. %
\begin{itemize}
\item[\upshape(NON)]%
  There are no O'Nan configurations in $\cH(C|R)$.
\item[\upshape(ALL)]%
  For any three distinct points $x,y,z$ on a block, there is a
  translation of $\cH(C|R)$ with center~$z$ mapping~$x$ to~$y$.
\end{itemize}
Let\/~$p$ be a point of\/ $\cH(C|R)$, and let $X,Y,Z$ be three
distinct blocks through~$p$. %
\begin{itemize}
\item[\upshape(TRA)] If $X,Y,Z$ intersect a block $B$ not through $p$, then
  for each point $z\ne p$ on either of the three blocks, there exists
  a (unique) block~$B'$ containing $z$ and intersecting the three
  blocks in three distinct points. %
  The block~$B'$ is the image of~$B$ under a translation of\/
  $\cH(C|R)$ with center~$p$, and each block through $p$ intersecting
  either of $B,B'$ intersects both $B$ and~$B'$.
\item[\upshape(TAN)] If $X,Y,Z$ intersect two disjoint blocks $B$ and $B'$ not
  containing $p$, then the intersection point of the two lines
  containing the blocks $B$ and $B'$ in the projective plane
  $\PG(2,C)$ is on the tangent line at~$p$.
\end{itemize}   
\end{lemma}

\begin{proof}
  Assertions (NON) and (ALL) have been proved
  in~\cite[Proposition\,2.7, Remark\,2.12]{MR4008656}. %
  For Assertion~(TRA), we apply the translation with center~$p$ that
  maps the intersection of the block joining~$p$ and~$z$ with~$B$ to
  the point~$z$.  %
  Uniqueness of~$B'$ is a consequence of (NON), %
  and the rest of Assertion~(TRA) follows from the fact that $B'$ is
  the image of~$B$ under a translation with center~$p$.
  
  It remains to verify Assertion~(TAN). %
  We have just seen that~$B'$ is the image of~$B$ under a translation
  of~$\cH(C|R)$ with center~$p$.  That translation is the restriction
  of an elation of~$\PG(2,C)$, the center is~$p$ and the axis is the
  tangent to~$\cH(C|R)$ at~$p$ in~$\PG(2,C)$. The intersection of the
  lines containing the blocks $B$ and $B'$ in the projective plane
  $\PG(2,C)$ is a point fixed by that elation, and thus contained in
  the tangent line at~$p$.
\end{proof}

\begin{defin}
  If $(X,\cB)$ and $(Y,\cD)$ are linear spaces then an embedding of
  $(X,\cB)$ into $(Y,\cD)$ is a pair of injective maps
  $\alpha\colon X\to Y$ and $\beta\colon \cB\to\cD$ such that $x\in X$
  and $B\in\cB$ are incident in $(X,\cB)$ exactly if $x^\alpha$ and
  $B^\beta$ are incident in $(Y,\cD)$. %
  We call $(X,\cB)$ a \emph{subunital} of~$(Y,\cD)$ if
  $(\operatorname{id},\operatorname{id})$ is an embedding of a unital
  $(X,\cB)$ into $(Y,\cD)$.
\end{defin}

\begin{defin}\label{def:standard}
  Let $E|F$ and $C|R$ be field extensions.  %
  An \emph{embedding} of field extensions from $E|F$ into $C|R$ is a
  field monomorphism $\eta\colon E\to C$ such that
  $F^\eta\subseteq R$. %
  An embedding $(\alpha,\beta)$ of $\cH(E|F)$ into $\cH(C|R)$ is
  called \emph{standard} if there is an embedding~$\eta$ of field
  extensions from $E|F$ into $C|R$ such that
  $\bigl((x,y,z)E\bigr)^\alpha = (x^\eta,y^\eta,z^\eta)C$ holds for
  each point $(x,y,z)E$ of~$\cH(E|F)$, up to some projective
  transformation of~$\PG(2,C)$.
\end{defin}

\section{Main Result}

\begin{lemma}\label{theoLem}
  Let $(P,\cL)$ be a linear space containing no O'Nan configuration,
  let $(U,\cB)$ be a unital of order~$q$, and let\/ $(\alpha,\beta)$
  be an embedding of\/ $(U,\cB)$ in $(P,\cL)$.  If two blocks $B_1$,
  $B_2$ of $(U,\cB)$ have no point of~$U$ in common, then $B_1^\beta$,
  $B_2^\beta$ are disjoint blocks of $(P,\cL)$.
\end{lemma}
\goodbreak%
\begin{proof}
  Aiming at a contradiction, suppose that two blocks $B_1,B_2\in\cB$
  are disjoint in $U$, but that $B_1^\beta$, $B_2^\beta$ contain a
  common point $x$. The absence of O'Nan configurations in $(P,\cL)$
  implies that two arbitrary blocks of $(U,\cB)$ both intersecting
  $B_1\cup B_2$ in exactly two points have no points off $B_1\cup B_2$
  in common. Hence the number of points in $U$ lying on a block
  intersecting $B_1\cup B_2$ in exactly two points is greater than
  $(q+1)^2(q-1) \geq q^3+1$; this is a contradiction. 
\end{proof}

\goodbreak%
\begin{thm}\label{sub}
  Let\/ $(U,\cB)$ be a finite subunital of the generalized hermitian
  unital\/ $\cH(C|R)$. Then $(U,\cB)$ is a hermitian unital, isomorphic
  to $\cH(E|F)$ for some quadratic extension $E|F$ of finite fields,
  and the embedding of\/ $(U,\cB) \cong \cH(E|F)$ into $\cH(C|R)$ is
  standard, %
  coming from an embedding of field extensions from $E|F$ into $C|R$.
  Moreover, the extension $C|R$ is separable, and $\cH(C|R)$ is a
  hermitian unital. %
\end{thm}

\begin{proof}
  Let $q$ denote the order of $(U,\cB)$. %
  In this proof, we suppose $q>2$; the case $q=2$ is treated
  separately in~\ref{noH2unusual} below.

  Let $p\in U$ be arbitrary, and let $B\in\cB$ be such that
  $p\notin B$. Let $B_0,B_1,\ldots,B_{q}$ be the blocks of $(U,\cB)$
  containing $p$ and intersecting $B$ nontrivially, say in
  $x_0,x_1,\ldots,x_q$, respectively. Let $x$ be an arbitrary point on
  $B_0\smallsetminus\{p,x_0\}$. We claim that at least one block of
  $(U,\cB)$ contains $x$ and intersects
  $B_1\cup B_2\cup\cdots\cup B_q$ in at least two points (different
  from~$p$).  Indeed, if not, then there are $q^2$ blocks through $x$
  different from $B_0$, which is a contradiction. So let $B_x$ be a
  block of $(U,\cB)$ containing at least three points (including~$x$)
  of $B_0\cup B_1\cup\cdots\cup B_q$. We note that $B_x$ and $B$ are
  disjoint by (NON). %
  For the same reason (or by~\ref{theoLem}) their extensions to
  $\cH(C|R)$ are also disjoint.  It then follows from~(TRA) that $B_x$
  intersects $B_i$ for each $i\in\{0,1,\ldots,q\}$, and~\ref{theoLem}
  yields that these intersection points belong to $U$. Hence we have
  shown that~(TRA) holds in the subunital~$(U,\cB)$.

  Now let $\tau$ be the translation of $\cH(C|R)$ with center $p$
  mapping $x_0$ to~$x$.  Let $y$ be any point of $U$ not
  on~$B_0$. Since $B$ was arbitrary, we may assume that $y\in B$, so
  without loss of generality $y=x_1$. By the uniqueness in~(TRA),
  $\tau$ maps $x_1$ to the intersection $B_x\cap B_1$. Since this
  intersection point belongs to $U$, it follows that $\tau$
  preserves $U$. Hence $(U,\cB)$ admits all translations and hence is
  hermitian by the main result of \cite{MR3090721}.

  Now consider the (standard) embedding of $\cH(C|R)$ in the
  projective plane $\PG(2,C)$. Then also $(U,\cB)$ is embedded in
  $\PG(2,C)$ and so by the Main Theorem of~\cite{MR4008656} there is a
  subfield $E\leq C$ of order $q^2$ and a subplane $\pi\cong\PG(2,E)$
  containing~$U$. %
  (Here we use $q>2$, for $q=2$ the Main Theorem of~\cite{MR4008656}
  does not apply.) %
  Hence there is a polarity $\rho_\pi$ of $\pi$ with absolute point
  set $U$. We now show that $\rho_\pi$ extends to a polarity $\rho$ of
  $\PG(2,C)$ where the absolute points are the points of $\cH(C|R)$.
  
  Consider the lines extending blocks of $(U,\cB)$ that meet at least
  three blocks through~$p$.  By~(TAN), any two of those lines have an
  intersection point on the tangent line~$T$ to~$U$ at~$p$
  in~$\pi$. Varying the blocks to be extended, we obtain more than one
  point on~$T$.  %
  The same description applies, \emph{mutatis mutandis}, to points on
  the tangent to $\cH(C|R)$ at~$p$ in~$\PG(2,C)$. So that line extends
  the tangent to $(U,\cB)$ in~$\pi$. %
  This already implies that not all tangent lines to $\cH(C|R)$
  contain the same point and so $C|R$ is separable
  by~\cite[2.3(2)]{MR4008656}. %
  Hence there is a polarity $\rho$ of $\PG(2,C)$ associated to
  $\cH(C|R)$. Since $U$ contains a quadrangle, and points of $U$ are
  mapped to lines of $\pi$ under the action of $\rho$, we see that
  $\rho$ preserves $\pi$. Since tangent lines to $(U,\cB)$ and
  $\cH(C|R)$ coincide in $\pi$, we see that $\rho_{|\pi} = \rho_\pi$.
  Hence the generator of the Galois group of $C|R$ preserves $E$ and
  induces $x\mapsto x^q$ in~$E$.

  This proves our main result completely for $q\neq 2$.  For $q=2$,
  see~\ref{noH2unusual}.
\end{proof}

If $C$ is finite of order $r^2$ then $E$ is unique with given order
$q^2$ and the field automorphism $x\mapsto x^r$ is not trivial on $E$,
which means that $E$ is not contained in the unique subfield $R$ of
order $r$; hence $C$ is an extension of $E$ of odd degree. %
Thus~\ref{sub} gives the following.

\begin{cor}
  If\/ $\UU$ is a unital of order~$q$ embedded in the hermitian unital
  $\cH(\FF_{r^2}|\FF_r)$ of order~$r$, then\/~$\UU$ is the hermitian
  unital~$\cH(\FF_{q^2}|\FF_q)$, there is an odd integer~$e$ such that
  $r=q^e$, and the embedding of the unital is standard.
\end{cor}

Up to projective equivalence, the embedded unital is obtained by using
a hermitian equation over $\FF_{q^2}|\FF_q$ to define the large unital
in $\PG(2,\FF_{r^2})$, and restricting coordinates to $\FF_{q^2}$ to
define the small unital.

\begin{rems}
  Wilbrink~\cite{MR690826} has characterized the finite hermitian
  unitals by three conditions. His condition~(I) is our~(NON)
  in~\ref{properties}.  Under the assumption~(NON), his condition~(II)
  is equivalent to our~(TRA). For unitals of even order, these two
  conditions alone suffice to characterize the hermitian unitals,
  see~\cite{MR3713359}.  Wilbrink's condition~(III) is too technical
  to state it here; compare~\cite[p.\,299]{MR4224387}.
\end{rems}

\section{The smallest unital}

\begin{thm}\label{noH2unusual}
  Let $C|R$ be a quadratic field extension. %
  The unital of order~$2$ is embedded in the generalized hermitian
  unital $\cH(C|R)$ if, and only if, $R$ has characteristic~$2$ and
  $C\cong R[X]/(X^2+X+1)$.  %
  The embedding is then a standard one.
\end{thm}
\begin{proof}
  Each unital of order~$2$ is isomorphic to the hermitian unital
  $\cH(\FF_4|\FF_2)$, and is clearly embedded in $\cH(C|R)$ if $R$ has
  characteristic~$2$ and $C\cong R[X]/(X^2+X+1)$.
  
  The unital of order~$2$ is also isomorphic to the affine plane of
  order~$3$ (see~\cite[10.16]{MR1189139}), and we know the embeddings
  of the latter into Moufang planes (\cite[3.7]{MR3856477},
  cp.~\cite[5.2]{MR901826}, \cite{MR2367750}): %
  Such an embedding is possible only if the coordinatizing alternative
  field contains an element $u$ with $u^2+u+1=0$; the points of
  the embedded plane are then given as
  $(0,0)$, $(1,0)$, $(0,1)$, $(-u,1)$, $(1,-u^2)$, $(-u,-u^2)$, $(0)$, 
  $(u)$, $(\infty)$, in suitable inhomogeneous
  coordinates.
  
  Consider a separable extension $C|R$ first, and assume that there is
  an embedding of the unital of order~$2$ into the hermitian unital
  $\cH(C|R)$. %
  We introduce homogeneous coordinates such that the points in
  question are \( [1,0,0]\), \([1,1,0]\), \([1,0,1]\), \([1,-u,1]\),
  \([1,1,-u^2]\), \([1,-u,-u^2]\), \([0,1,0]\), \([0,1,u]\), and
  \([0,0,1]\). %
  In those homogeneous coordinates, the hermitian unital consists of
  all points $[x,y,z]$ satisfying the equation
  $(\bar{x},\bar{y},\bar{z})M(x,y,z)^T=0$, where $x\mapsto\bar{x}$ is
  the generator of the Galois group of $C|R$, and $M$ is a
  non-singular hermitian $3\times3$ matrix over~$C$. %
  (Multiplying $M$ by a skew-symmetric field element swaps hermitian with
  skew-hermitian forms but preserves the unital, so we are in
  accordance with~\ref{identify}.) %
  Now $[1,0,0],[0,1,0],[0,0,1]\in\cH(C|R)$ implies that each diagonal
  entry of~$M$ is zero. Thus $M=\left(
    \begin{smallmatrix}
      0 & a & \bar{b} \\
      \bar{a} & 0 & c \\
      b & \bar{c} & 0
    \end{smallmatrix}
  \right)$ and the equation becomes
  \[ %
    \operatorname{tr}(\bar{x}ay) + \operatorname{tr}(\bar{z}bx) +
    \operatorname{tr}(\bar{y}cz) =0 \,;
  \]
  where $\operatorname{tr}(w):= w+\bar{w}$ is the trace of~$w$.
  
  Evaluating this equation for \([1,1,0]\), \([1,0,1]\), and
  \([0,1,u]\), respectively, we
  obtain $\operatorname{tr}(a) = \operatorname{tr}(b) =
  \operatorname{tr}(cu) =0$.

  If $\bar{u}=u$ then $0=\operatorname{tr}(cu) = (\bar{c}+c)u$ yields
  $\operatorname{tr}(c)=0$ and then $\det{M} = \operatorname{tr}(abc)
  = ab\operatorname{tr}(c) = 0$, a
  contradiction. %
  So $\bar{u}\ne u$ is another root of $X^2+X+1$. This entails
  $\bar{u}=-u-1=u^2$ and $\bar{u}^2=u$. %

  Now we evaluate the equation above
  for \([1,1,-u^2]=[1,1,-\bar{u}]\), \([1,-u,1]\), and
  \([1,-u,-u^2]=[1,-u,-\bar{u}]\), respectively, and obtain
  \[
    \begin{array}{rcccl}
      0 &=& \operatorname{tr}(a)-\operatorname{tr}(\bar{u}b)-\operatorname{tr}(\bar{u}c) %
      &=&  b(\bar{u}-u) - \operatorname{tr}(\bar{u}c)  \\
      0 &=& -\operatorname{tr}(au)+ \operatorname{tr}(b)-\operatorname{tr}(\bar{u}c) %
      &=& a(\bar{u}-u) - \operatorname{tr}(\bar{u}c) \\
      0 &=& -\operatorname{tr}(au)-\operatorname{tr}(ub)+\operatorname{tr}(c\bar{u}^2) %
      &=& (a+b)(\bar{u}-u)          \,.
    \end{array}
  \]
  As $\bar{u}-u\ne0$, these equations give $a=b$ and then $0=2a$. %
  So $R$ has characteristic~$2$, and $\FF_2(u) \cong\FF_4$ is
  contained in $C$ but not in~$R$. %
  We obtain an embedding of the field extension $\FF_{4}|\FF_2$ into
  $C|R$, and the embedding of the unital~$\cH(\FF_4|\FF_2)$ of order 2
  in $\cH(C|R)$ is a standard embedding; see~\ref{def:standard}.

  It remains to show that the unital of order~$2$ is not embedded in
  $\cH(C|R)$ if $C|R$ is inseparable. %
  As above, we know from~\cite[3.7]{MR3856477} that any embedding of
  $\cH(\FF_4|\FF_2)$ into~$\PG(2,C)$ is projectively equivalent to the
  one where the embedded points have coordinates in $\{0,1,u,u^2\}$,
  with $u\in C$ and $u^2=u+1$; we use that $C$ has characteristic~$2$
  because $C|R$ is inseparable. Now the set $\{0,1,u,u^2\}$ is a
  subfield of order~$4$ in~$C$, and it is contained in~$R$ because
  otherwise the extension $C|R$ would be separable. Thus the embedded
  unital of order~$2$ is contained in a finite subplane
  $\pi\cong\PG(2,\FF_4)$ of~$\PG(2,C)$.

  Let $B$ be a block of~$\cH(C|R)$ joining two points of the embedded
  unital of order~$2$. Then $B\cap\pi$ is a set of three points. We
  introduce coordinates for the line~$L$ of~$\PG(2,C)$ containing~$B$
  in such a way that $0$, $1$, and~$\infty$ are the coordinates of the
  three points in~$B\cap\pi$. %
  Then~$B$ consists of the points with coordinates in
  $R\cup\{\infty\}$ because the blocks of~$\cH(C|R)$ are Baer sublines
  with respect to $C|R$, see~\cite[2.8]{MR4008656}. Thus $B\cap\pi$
  contains all points of~$L\cap\pi$ (i.e., those with coordinates $0$,
  $1$, $u$, $u^2$, and~$\infty$), contradicting the fact that
  $|B\cap\pi|=3$. This contradiction yields that there is no embedding
  of $\cH(\FF_4|\FF_2)$ in $\cH(C|R)$ if $C|R$ is inseparable.
\end{proof}

\goodbreak%
%%%%%%%%%%%%%%%%%%%%%%%%%%%%%%%%%%%%%%%%%%%%%%%%%%%%%%%%%%%%%%%%%%%%%%%
% \bibliographystyle{mybibstyle-noISSN-noISBN} %%% %%%%% %%%%% %%%%% %%%%
% \bibliography{myBibliography}  %%%%% %%%%% %%%%% %%%%% %%%%% %%%%% %%%%

\begin{thebibliography}{10}
\providecommand{\href}[2]{#2}
\providecommand{\eprint}[1]{\href{http://arxiv.org/abs/#1}{#1}}
\providecommand{\bbldiplomarbeit}{Diplomarbeit}
\def\bbland{and}                \def\bbletal{et~al.}
\def\bbleditors{editors}        \def\bbleds{eds.}
\def\bbleditor{editor}          \def\bbled{ed.}
\def\bbledby{edited by}
\def\bbledition{edition}        \def\bbledn{edn.}
\def\bblvolume{volume}          \def\bblvol{vol.}
\def\bblof{of}
\def\bblnumber{number}          \def\bblno{no.}
\def\bblin{in}
\def\bblpages{pages}            \def\bblpp{pp.}
\def\bblpage{page}              \def\bblp{p.}
\def\bbleidpp{pages}
\def\bblchapter{chapter}        \def\bblchap{chap.}
\def\bbltechreport{Technical Report}
\def\bbltechrep{Tech. Rep.}
\def\bblmthesis{Master's thesis}
\def\bblphdthesis{Ph.D. thesis}
\def\bblfirst{First}            \def\bblfirsto{1st}
\def\bblsecond{Second}          \def\bblsecondo{2nd}
\def\bblthird{Third}            \def\bblthirdo{3rd}
\def\bblfourth{Fourth}          \def\bblfourtho{4th}
\def\bblfifth{Fifth}            \def\bblfiftho{5th}
\def\bblst{st}  \def\bblnd{nd}  \def\bblrd{rd}
\def\bblth{th}
\def\bbljan{January}  \def\bblfeb{February}  \def\bblmar{March}
\def\bblapr{April}    \def\bblmay{May}       \def\bbljun{June}
\def\bbljul{July}     \def\bblaug{August}    \def\bblsep{September}
\def\bbloct{October}  \def\bblnov{November}  \def\bbldec{December}
\newcommand{\Capitalize}[1]{\uppercase{#1}}
\newcommand{\capitalize}[1]{\expandafter\Capitalize#1}

%% definitions for hyperlinks used in mybibstyle.bst etc.
\providecommand{\url}[1]{\href{#1}{#1}}
\providecommand{\urlprefix}{}
\let\oldunderscore_
\catcode`\_=13
\providecommand{\doi}[1]{\href{http://dx.doi.org/#1}{{\def_{\_}\normalfont\ttfamily doi:#1}}\let_\oldunderscore}
% we adapt \MR , making it work well with MathSciNet bibtex version:
\providecommand{\MR}[1]{\relax\ifhmode\unskip\space\fi \MRnumberextract#1 \,}
\def\MRnumberextract#1 #2\,{\MRhref{#1}{#2}}%
% \MRhref is called by the amsart/book/proc definition of \MR.
\providecommand{\MRhref}[2]{%
  \href{http://www.ams.org/mathscinet-getitem?mr=#1}{MR\,#1 #2}}
\providecommand{\ZBL}[1]{\relax\ifhmode\unskip\space\fi \ZBLhref{#1}}
\providecommand{\ZBLhref}[1]{%
  \href{http://zbmath.org/?q=an:#1}{Zbl #1}}
\providecommand{\JfM}[1]{\relax\ifhmode\unskip\space\fi \JfMhref{#1}}
\providecommand{\JfMhref}[1]{%
  % \href{https://www.emis.de/cgi-bin/jfmen/MATH/JFM/quick.html?type=html&an=#1}{JfM #1}
  \href{http://zbmath.org/?q=an:#1}{JfM #1}}

\bibitem{MR901826}
M.~S. Abdul-Elah, M.~W. Al-Dhahir, \bbland{} D.~Jungnickel, \emph{{$8_3$} in
  {${\rm PG}(2,q)$}}, Arch. Math. (Basel) \textbf{49} (1987), \bblno{}~2,
  141--150, \doi{10.1007/BF01200478}. \MR{901826.} \ZBL{0595.51003}.

\bibitem{MR2440325}
S.~G. Barwick \bbland{} G.~Ebert, \emph{Unitals in projective planes}, Springer
  Monographs in Mathematics, Springer, New York, 2008,
  \doi{10.1007/978-0-387-76366-8}. \MR{2440325.} \ZBL{1156.51006}.

\bibitem{MR3090721}
T.~Grundh{\"o}fer, M.~J. Stroppel, \bbland{} H.~{Van Maldeghem}, \emph{Unitals
  admitting all translations}, J. Combin. Des. \textbf{21} (2013), \bblno{}~10,
  419--431, \doi{10.1002/jcd.21329}. \MR{3090721.} \ZBL{1276.05021}.

\bibitem{MR4008656}
T.~Grundh{\"o}fer, M.~J. Stroppel, \bbland{} H.~{Van Maldeghem},
  \emph{Embeddings of hermitian unitals into pappian projective planes},
  Aequationes Math. \textbf{93} (2019), \bblno{}~5, 927--953,
  \doi{10.1007/s00010-019-00652-x}. \MR{4008656.} \ZBL{07109887}.

\bibitem{MR4224387}
T.~Grundh{\"o}fer, M.~J. Stroppel, \bbland{} H.~{Van Maldeghem},
  \emph{Embeddings of unitals such that each block is a subline}, Australas. J.
  Combin. \textbf{79} (2021), 295--301,
  \urlprefix\url{https://ajc.maths.uq.edu.au/pdf/79/ajc_v79_p295.pdf}.
  \MR{4224387.} \ZBL{7352033}.

\bibitem{MR0333959}
D.~R. Hughes \bbland{} F.~C. Piper, \emph{Projective planes}, {G}raduate
  {T}exts in {M}athematics ~6, Springer-Verlag, New York, 1973. \MR{0333959.}
  \ZBL{0484.51011}.

\bibitem{MR3713359}
A.~M.~W. Hui, \emph{A geometric proof of {W}ilbrink's characterization of even
  order classical unitals}, Innov. Incidence Geom. \textbf{15} (2017),
  145--167, \doi{10.2140/iig.2017.15.145}. \MR{3713359.} \ZBL{1415.51009}.

\bibitem{MR0295934}
M.~E. O'Nan, \emph{Automorphisms of unitary block designs}, J. Algebra
  \textbf{20} (1972), 495--511, \doi{10.1016/0021-8693(72)90070-1}.
  \MR{0295934.} \ZBL{0241.05013}.

\bibitem{MR2367750}
G.~Pickert, \emph{Near-embeddings of the affine plane with 9 points into
  {D}esarguesian projective and affine planes}, Note Mat. \textbf{27} (2007),
  \bblno{}~1, 11--19, \doi{10.1285/i15900932v27n1p11}. \MR{2367750.}
  \ZBL{1195.51009}.

\bibitem{MR3856477}
M.~J. Stroppel, \emph{Generalizing the {P}appus and {R}eye configurations},
  Australas. J. Combin. \textbf{72} (2018), 249--272,
  \urlprefix\url{http://ajc.maths.uq.edu.au/pdf/72/ajc_v72_p249.pdf}.
  \MR{3856477.} \ZBL{07021381}.

\bibitem{MR1189139}
D.~E. Taylor, \emph{The geometry of the classical groups}, Sigma Series in Pure
  Mathematics ~9, Heldermann Verlag, Berlin, 1992. \MR{1189139.}
  \ZBL{0767.20001}.

\bibitem{MR1725957}
H.~{Van Maldeghem}, \emph{Generalized polygons}, Monographs in Mathematics ~93,
  Birkh\"auser Verlag, Basel, 1998, \doi{10.1007/978-3-0348-0271-0}.
  \MR{1725957.} \ZBL{0914.51005}.

\bibitem{MR690826}
H.~A. Wilbrink, \emph{A characterization of the classical unitals}, \bblin{}
  \emph{Finite geometries ({P}ullman, {W}ash., 1981)}, Lecture Notes in Pure
  and Appl. Math. ~82, \bblpp{} 445--454, Dekker, New York, 1983.
  \urlprefix\url{https://ir.cwi.nl/pub/6786}. \MR{690826.} \ZBL{0507.05014}.

\end{thebibliography}
% \end{document}
%%%%%%%%%%%%%%%%%%%%%%%%%%%%%%%%%%%%%%%%%%%%%%%%%%%%%%%%%%%%%%%%%%%%%%% 
\providecommand{\noopsort}[1]{}\def\cprime{$'$}
  \def\polhk#1{\setbox0=\hbox{#1}{\ooalign{\hidewidth
  \lower1.5ex\hbox{`}\hidewidth\crcr\unhbox0}}}

%%%%%%%%%%%%%%%%%%%%%%%%%%%%%%%%%%%%%%%%%%%%%%%%%%%%%%%%%%%%%%%%%%%%%%% 
\end{document}